\documentclass[11pt]{amsart}

\usepackage{amssymb}
\usepackage{latexsym}
\usepackage{amsmath}

\begin{document}

\newtheorem{theorem}{Theorem}[section]
\newtheorem{proposition}{Proposition}[section]
\newtheorem{lemma}{Lemma}[section]
\newtheorem{definition}{Definition}[section]
\newtheorem{problem}{Problem}[section]

\title{Spectral measures of small index principal graphs}
\author{Teodor Banica}
\address{T.B.: Department of Mathematics, Universit{\'e} Paul Sabatier, 
118 route de Narbonne, 31062 Toulouse, France}
\email{banica@picard.ups-tlse.fr}
\author{Dietmar Bisch}
\thanks{D.B. was supported by NSF under Grant No. DMS-0301173}
\address{D.B.: Department of Mathematics, Vanderbilt University, 
1326 Stevenson Center, Nashville, TN 37240, USA}
\email{dietmar.bisch@vanderbilt.edu}
\subjclass[2000]{46L37 (46L54)}
\keywords{Subfactor, Principal graph, ADE classification, Semicircle law}

\begin{abstract}
The principal graph $X$ of a subfactor with finite Jones index
is one of the important algebraic invariants of the subfactor.
If $\Delta$ is the adjacency matrix of $X$ we consider the equation 
$\Delta=U+U^{-1}$. When $X$ has square norm $\leq 4$ the spectral 
measure of $U$ can be averaged by using the map $u\to u^{-1}$, and 
we get a probability measure $\varepsilon$ on the unit circle which 
does not depend on $U$. We find explicit formulae for this measure 
$\varepsilon$ for the principal
graphs of subfactors with index $\le 4$, the (extended)
Coxeter-Dynkin graphs of type $A$, $D$ and $E$. The moment generating
function of $\varepsilon$ is closely related to Jones' $\Theta$-series.
\end{abstract}
\maketitle

\section*{Introduction}

The Coxeter-Dynkin graphs of type A, D$_{\text{\rm even}}$, E$_6$,
E$_8$, and the extended Coxeter-Dynkin graphs of type ADE
appear in the theory of 
subfactor as basic invariants for inclusions of II$_1$ factors 
with Jones index $\le 4$ (\cite{jo1}, \cite{oc}, \cite{po1}, \cite{po2},
see also \cite{ek}, \cite{ghj}). They are fusion graphs 
of subfactor representations
and capture the algebraic information contained in the standard
invariant of the subfactor $N \subset M$ (see e.g. \cite{bi}, 
\cite{po2}, \cite{ek}, \cite{ghj}). Such a graph $X$ is bipartite and has 
a distinguished vertex $1$. Of particular interest in subfactor 
theory are the number of length $2k$ loops on $X$ based at $1$, 
since these are the dimensions of the higher relative commutants
associated to $N \subset M$. If combined in a formal power series
$f(z)$ we obtain the Poincar\'e series of $X$. 

A related series $\Theta(u)$, with  $z^{1/2}=u+u^{-1}$, is considered 
by Jones in \cite{jo3}. Jones made the remarkable discovery that if
the subfactor has index $> 4$ then the coefficients of this series are 
necessarily {\it positive} integers. 
The series $f$ and $\Theta$ are natural invariants of the subfactor.
In this paper we compute explicitly measures whose generating
series of moments are (essentially) these two power series in the case
the graphs have square norm $\le 4$. These
measures can then be regarded as invariants of the subfactors.

For a graph $X$ of square norm $\leq 4$ we have the following 
measure-theoretic version of $\Theta$. Consider the equation 
$\Delta=U+U^{-1}$, 
where $\Delta$ is the adjacency matrix of $X$. We can 
average the spectral measure of $U$ by using the map $u\to u^{-1}$, 
and we get a probability measure 
$\varepsilon$ on the unit circle which does not depend on $U$.

We find that for the (extended) Coxeter-Dynkin graphs of type 
A and D this measure $\varepsilon$ is given by very simple formulae
as follows:
\begin{eqnarray*}
A_{n-1}&\rightarrow&\alpha\,d_n\cr
D_{n+1}&\rightarrow&\alpha\,d^\prime_n\cr
A_\infty&\rightarrow&\alpha\,d\cr
A_{2n}^{(1)}&\rightarrow&d_n\cr
A_{-\infty,\infty}&\rightarrow&d\cr
D_{n+2}^{(1)}&\rightarrow&(d^\prime_1+d_n)/2\cr
D_\infty&\rightarrow&(d^\prime_1+d)/2
\end{eqnarray*}

Here $d$, $d_n$, $d^\prime_n$ are the uniform measures on the unit circle,
the uniform measure supported on the $2n$-th roots of unity, and the 
uniform measure supported on the $4n$-th roots of unity of odd order. 
The fundamental density $\alpha$ is given by 
$\alpha(u)=2Im(u)^2$ and corresponds via $x=u+u^{-1}$ to the semicircle 
law from \cite{vdn}.

For the graphs of type $E_6$, $E_7$, $E_8$, $E_6^{(1)}$, $E_7^{(1)}$, 
and $E_8^{(1)}$, $\varepsilon$ is given by the following formulae:
\begin{eqnarray*}
E_6&\to&\alpha\,d_{12}+(d_{12}-d_6-d_4+d_3)/2\cr
E_7&\to&\varepsilon_7\cr
E_8&\to&\varepsilon_8\cr
E_6^{(1)}&\to&\alpha\,d_3+(d_2-d_3)/2\cr
E_7^{(1)}&\to&\alpha\,d_4+(d_3-d_4)/2\cr
E_8^{(1)}&\to&\alpha\,d_6+(d_5-d_6)/2
\end{eqnarray*}

Here $\varepsilon_7$, $\varepsilon_8$ denote certain exceptional 
measures for which we do not have closed formulae (we do compute
their moment generating series though). These results are obtained 
by explicit computations.
They could also be obtained by using planar algebras methods, see 
\cite{jo3}, \cite{jo-re}, \cite{re} for related work.

The measure $\varepsilon$ computed in this article should be viewed as 
an analytic invariant for subfactors of index $\leq 4$. It is unclear 
what the appropriate generalization of $\varepsilon$ to subfactors 
of index $> 4$ should be. However, in light of Jones' work in 
\cite{jo3}, the measure $\varepsilon$ should be related to certain 
representations of planar algebras. It should shed some light on the 
structure of subfactor planar algebras (or equivalently standard 
invariants) arising from subfactors with index $>4$. We intend to
come back to this question in future work.

Similar considerations make sense for quantum groups. The hope 
would be that an analytic invariant for quantum groups would emerge 
from the Weingarten formula in \cite{bc}, and work here is in progress. 
This in turn might be related to the results in the present article 
via Di Francesco's formula in \cite{df}.

The paper is dedicated to the proofs of the above results and is organized 
as follows. Section 1 fixes the notation and contains some preliminaries.
In sections 2, 3 and 4 we divide the graphs of type A and D into three 
classes -- circulant graphs, graphs with A$_n$ tails, and graphs with 
fork tails -- and we compute $\varepsilon$. In sections 5, 6 and 7 we 
discuss the graphs of type E using a key lemma in section 3.

\subsection*{Acknowledgements}
Most of this paper was written while T.B.
was visiting Vanderbilt University and D.B. was visiting the 
Universit{\'e} Paul Sabatier. The authors are grateful for the 
hospitality received at both institutions.

\section{Spectral measures and the Jones series}

We collect in this section several known results about spectral measures
associated to graphs and their Stieltjes transforms. We relate them to
a natural power series discovered by Jones in the context of planar
algebras and their representation theory (\cite{jo2}, \cite{jo3}).

Let $X$ be a (possibly infinite) bipartite graph with 
distinguished vertex labeled by "$1$". Since $X$ 
is bipartite, its adjacency matrix $\Delta$ is given by
$$\Delta=\begin{pmatrix}0&M\cr M^t&0\end{pmatrix}$$
where $M$ is a rectangular matrix with non-negative 
integer entries. If we let $L=MM^t$ and $N=M^tM$, then
$$\Delta^2=\begin{pmatrix}L&0\cr 0&N\end{pmatrix}$$

For a matrix $T$ with entries labeled by the vertices of $X$, we use the 
following notation ("1" is the label of the distinguished vertex):
$$\int T=T_{11}$$

We call $T_{11}$ the integral of the matrix $T$. 

\begin{definition}
The spectral measure of $X$ is the probability measure 
$\mu$ on ${\mathbb R}$ satisfying
$$\int_{-\infty}^\infty \varphi (x)\,d\mu (x)=\int\varphi(\Delta)$$
for any continuous function $\varphi :{\mathbb R}\to{\mathbb C}$.
\end{definition}

Note that the spectral measure of $\Delta$ can be regarded as an 
invariant of $X$.
The spectral measure is uniquely determined by its moments. The generating 
series of these moments is called the Stieltjes transform 
$\sigma$ of $\mu$, i.e.
$$\sigma(z)=\int_{-\infty}^\infty\frac{1}{1-zx}\,d\mu(x).$$

This is related to the Poincar\'e series of $X$, which appears as 
the generating function of the numbers ${\rm loop}(2k)$, counting 
loops of length $2k$ on $X$ based at $1$. The Poincar\'e series
of $X$ is defined as

$$f(z)=\sum_{k=0}^\infty \, {\rm loop}(2k) z^k.$$

We have then the following well-known result. 

\begin{proposition} Let $X$ be a bipartite graph with spectral 
measure $\mu$ and Poincar\'e series $f$.
The  Stieltjes transform $\sigma$ of $\mu$ is given by
$\sigma(z)=f(z^2)$.
\end{proposition}

\begin{proof}
We compute $\sigma(z)$ using the fact that the integral of 
$\Delta^l$ is ${\rm loop}(l)$.
\begin{eqnarray*}
\sigma(z)
=\int\frac{1}{1-z\Delta}
=\sum_{l=0}^\infty {\rm loop}(l)z^l
\end{eqnarray*}

Since ${\rm loop}(l)=0$ for $l$ odd, we have
$$f(z^2)=\sum_{k=0}^\infty {\rm loop}(2k)z^{2k}=\sum_{l=0}^\infty 
{\rm loop}(l)z^l.$$

This proves the statement.
\end{proof}

We assume now that the graph $X$ has norm $\leq 2$, that is the
matrix $\Delta$ has norm $\leq 2$. Thus the support of $\mu$, 
which is contained in the 
spectrum of $\Delta$, is contained in $[-2,2]$.

Let ${\mathbb T}$ be the unit circle, and consider the following map 
$\Phi :{\mathbb T}\to [-2,2]$, defined by $\Phi (u)=u+u^{-1}$.
Any probability measure $\varepsilon$ on ${\mathbb T}$ produces a 
probability measure $\mu =\Phi_*(\varepsilon)$ on $[-2,2]$, 
according to the following formula, valid for any continuous 
function $\varphi :[-2,2]\to{\mathbb C}$:

$$\int_{\mathbb R}\varphi(x)\,d\mu(x)=\int_{\mathbb T}\varphi(u+u^{-1})
\, d\varepsilon(u)$$

\medskip
We can obtain in this way all probability measures on $[-2,2]$. Given $\mu$, 
there is a unique probability measure $\varepsilon$ on $\mathbb T$ 
satisfying $\Phi_*(\varepsilon)=\mu$ with the normalisation 
$d\varepsilon (u)=d\varepsilon (u^{-1})$. 

\begin{definition}
The spectral measure of $X$ (on ${\mathbb T}$) is the 
probability measure 
$\varepsilon$ on ${\mathbb T}$ given by
$$\int_{\mathbb T} \varphi (u+u^{-1})\,d\varepsilon (u)=\int\varphi(\Delta)$$
for any continuous function $\varphi :[-2,2]\to{\mathbb C}$, with the 
normalisation $d\varepsilon (u)=d\varepsilon (u^{-1})$.
\end{definition}

The generating series of the moments of $\varepsilon$ (the Stieltjes 
transform) is given by 
$$S(q)=\int_{\mathbb{T}} \frac{1}{1-qu}\,d\varepsilon (u).$$

\medskip
Following Jones (\cite{jo2}, \cite{jo3}), given a subfactor planar 
algebra ($P =$ ($P_0^\pm$, ($P_k$)$_{k \ge 1}$) with parameter
$\delta$, the associated Poincar{\'e} series is defined as 

$$f(z) = \frac{1}{2}(\dim P_0^+ +\dim P_0^-)+ 
\sum_{k=1}^\infty \dim P_k z^k.$$ 

Jones introduced in \cite{jo3} an associated series $\Theta$, which
is essentially obtained from the Poincar\'e series by a change of 
variables $ z \to \frac{q}{1+q^2}$. In this paper we will call this 
series the {\it Jones series}. If $\delta > 2$, then the Jones series 
is the dimension generating function for the multiplicities of 
certain Temperley-Lieb modules which appear in the decomposition of 
the planar algebra $P$ 
viewed as a module for the Temperley-Lieb planar algebra (\cite{jo3},
see also \cite{jo-re}, \cite{re}).
Using the formula for the $\Theta$-series from \cite{jo3}
we define the {\it Jones series of the graph} $X$ by 

$$\Theta (q)=q+\frac{1-q}{1+q}\, f\left( \frac{q}{(1+q)^2}\right)$$
where $f(z)$ is the Poincar\'e series of $X$.
With these notations, we have then the following result.

\begin{proposition} 
Let $X$ be a bipartite graph with spectral
measure $\varepsilon$ (on $\mathbb{T}$) and Jones series $\Theta$.
The  Stieltjes transform of $\varepsilon$ is given by
$2S(q)=\Theta(q^2)-q^2+1$.
\end{proposition}

\begin{proof}
We compute $S$ in terms of $\varepsilon$.
\begin{eqnarray*}
2S(q)
&=&\int_{\mathbb{T}}\frac{1}{1-qu}\,d\varepsilon (u)+
\int_{\mathbb{T}}\frac{1}{1-qu^{-1}}\,d\varepsilon (u)\cr
&=&1+\int_{\mathbb{T}}\frac{1-q^2}{1-q(u+u^{-1})+q^2}\,d\varepsilon (u)
\end{eqnarray*}

We compute now $\Theta$ in terms of $\mu$.
\begin{eqnarray*}
\Theta(q^2)-q^2
&=&\frac{1-q^2}{1+q^2}\, f\left( \frac{q^2}{(1+q^2)^2}\right)\cr
&=&\int_{-\infty}^\infty\frac{1-q^2}{1-qx+q^2}\,d\mu(x)
\end{eqnarray*}

This formula in the statement follows now from the definition of $\varepsilon$.
\end{proof}

In the next sections we compute explicitly the spectral measure
(on $\mathbb{T}$) for graphs which appear in the classification
of subfactors with index $\le 4$ (the so-called {\it principal
graphs}). We use standard notation for 
these graphs (see for instance \cite{ghj}).

\section{Circulant graphs}

Consider the circulant graph $A_{2n}^{(1)}$, that is
$A_{2n}^{(1)}$ is the $2n$-gon, and choose any vertex as the distinguished 
vertex $1$.

\begin{theorem}
Let $X$ be the circulant graph $A_{2n}^{(1)}$. The spectral measure 
of $X$ (on ${\mathbb T}$) is given by
$$d\varepsilon (u)=d_nu$$
where $d_n$ is the uniform measure on $2n$-th roots of unity.
\end{theorem}

\begin{proof}
We identify $X$ with the group $\{ w^k\}_{0\le k \le 2n-1}$ of 
$2n$-th roots of unity, 
where $w=e^{i\pi/{n}}$. The adjacency matrix of $X$ acts on functions 
$f\in{\mathbb C}(X)$ in the following way.
$$\Delta f(w^s)=f(w^{s-1})+f(w^{s+1})$$

Consider the following operators $U$ and $U^{-1}$:
$$Uf(w^s)=f(w^{s+1})$$
$$U^{-1}f(w^s)=f(w^{s-1})$$

We have $\Delta=U+U^{-1}$. The moments of the spectral measure of $U$ are
obtained as follows:
$$\int U^k=<U^k\delta_1,\delta_1>=<\delta_{w^k},\delta_1>=\delta_{w^k,1}$$

We compute now the moments of the measure in the statement.
$$\int_{\mathbb T} u^k\,d\varepsilon(u)=
\int_{\mathbb T} u^k\,d_nu=\delta_{w^k,1}$$

Thus $\varepsilon$ is the spectral measure of $U$, and together with 
the identity $d\varepsilon(u)=d\varepsilon(u^{-1})$, we get the result.
\end{proof}

The graph $A_{-\infty,\infty}$ is the set ${\mathbb Z}$ with consecutive 
integers connected by edges. Choose any vertex as the distinguished vertex 
labeled $1$. 

\begin{theorem}
Let $X$ be the graph $A_{-\infty,\infty}$. The spectral measure 
of $X$ (on ${\mathbb T}$) is given by
$$d\varepsilon (u)=du$$
where $du$ is the uniform measure on the unit circle.
\end{theorem}

\begin{proof}
This follows from Theorem 2.1 by letting $n\to\infty$, or from a direct 
loop count.
\end{proof}

\section{Graphs with tails}

The Coxeter-Dynkin graph of type A$_n$, $n \ge 2$, has $n$ vertices
and the distinguished vertex $1$ labels a vertex at one end of the graph, 
i.e. A$_n$ is bipartite graph of the form

$$1-\alpha_1 -2-\alpha_2 -3 - \alpha_3 :::\beta$$
where $1$, $2$, $3$, $\dots$ labels the even vertices and
$\alpha_1$, $\alpha_2$, $\dots$ labels the odd vertices. $\beta$
is either even or odd depending on the parity of $n$.

Consider a sequence of graphs $X_k$ obtained by adding $A_{2k}$ tails 
to a finite graph $\Gamma$. We let
$$X_0=1:::\Gamma$$
$$X_1=1-\alpha -2:::\Gamma$$
$$X_2=1-\alpha -2-\beta -3:::\Gamma$$
where $1$ denotes the distinguished vertex, $\alpha$, $2$, $\beta$ and
$3$, $\dots$ denote vertices connected by single edges as indicated and
$\Gamma$ denotes a finite graph connected by a single edge to the
preceding vertex. For instance, $X_1$ is obtained by attaching
$A_3$ to $\Gamma$ (since the vertex 2 of $A_3$ is identified with
one of the vertices of $\Gamma$, we have attached an $A_2$-tail to
$\Gamma$), $X_2$ by attaching $A_5$ to $\Gamma$ etc. Similarly we
define the graph $X_k$.

We denote by $L_0$ the matrix appearing on the top left of the square of 
the adjacency matrix of $X_0$, that is:
$$\Delta_0^2=\begin{pmatrix}L_0&0\cr 0&N_0\end{pmatrix}$$

We compute the Jones series of each $X_k$ in the next lemma.

\begin{lemma}
The Jones series of the graphs $X_{k}$, $k \ge 0$, is given by
$$\frac{\Theta(q)-q}{1-q}=\frac{1-Pq^{2k}}{1-Pq^{2k+1}}$$
where $P$ is defined by the formula
$$P=\frac{P_1-q^{-1}P_0}{P_1-qP_0}$$
and $P_i=P_i(y)=\det (y-K_i)$, $i=0$, $1$, with $y=2+q+q^{-1}$, 
and with $K_0$, $K_1$ being the following matrices:

-- $K_0$ is obtained from $L_0$ by deleting the first row and column.

-- $K_1$ is obtained from $L_0$ by adding $1$ to the first entry.
\end{lemma}

\begin{proof} We use the notation fixed in section 1.
The matrix $M_k$ in the adjacency matrix of $X_k$ is given by a 
row vector $w$ and a matrix $M$, as follows:
$$M_0=\begin{pmatrix}w\cr M\end{pmatrix}\hskip 5mm
M_1=\begin{pmatrix}1&\cr1&w\cr&M\end{pmatrix}\hskip 5mm
M_2=\begin{pmatrix}1&&\cr1&1&\cr&1&w\cr&&M\end{pmatrix}$$

The corresponding matrices $L_k$ are given by the real number $a=ww^t+1$, 
the row vector $u=wM^t$, the column vector $v=Mw^t$, and the matrix 
$N=MM^t$ as follows:
$$L_0=\begin{pmatrix}w\cr M\end{pmatrix}\begin{pmatrix}w^t&M^t\end{pmatrix}=\begin{pmatrix}a-1&u\cr v&N\end{pmatrix}$$
$$L_1=
\begin{pmatrix}1&\cr1&w\cr&M\end{pmatrix}
\begin{pmatrix}1&1&\cr&w^t&M^t\end{pmatrix}=
\begin{pmatrix}1&1&\cr1&a&u\cr&v&N\end{pmatrix}$$
$$L_2=
\begin{pmatrix}1&&\cr1&1&\cr&1&w\cr&&M\end{pmatrix}
\begin{pmatrix}1&1&&\cr&1&1&\cr&&w^t&M^t\end{pmatrix}=
\begin{pmatrix}1&1&&\cr1&2&1&\cr&1&a&u\cr&&v&N\end{pmatrix}$$

It is now clear what the form of the matrices $M_k$ and $L_k$ is for 
general $k$. Consider the matrix $K_k$, $k \ge 0$, obtained from $L_k$ 
by deleting the first row and column, in other words
$$K_0=(N)\hskip 5mm K_1=\begin{pmatrix}a&u\cr v&N\end{pmatrix}\hskip 5mm
K_2=\begin{pmatrix}2&1&\cr1&a&u\cr&v&N\end{pmatrix}$$
and similarly for general $k$.

We have the following formula for the Poincar\'e series of $X_k$, 
with $y=z^{-1}$:

\begin{eqnarray*}
f_k(z)
=\int\frac{1}{1-zL_k}
=\frac{\det (1-zK_k)}{\det (1-zL_{k})}
=y\cdot\frac{\det (y-K_{k})}{\det (y-L_{k})}\cr
\end{eqnarray*}

The characteristic polynomials $P_k=\det (y-K_k)$ and $Q_k=\det (y-L_k)$ 
satisfy the following two identities, obtained by developing determinants 
at the top left.

$$P_{k+1}=(y-2)P_{k}-P_{k-1}$$
$$Q_{k}=(y-1)P_{k}-P_{k-1}$$

We consider the first identity and the second identity minus the first, 
with the change of variables $y=2+q+q^{-1}$, and obtain
$$P_{k+1}=(q+q^{-1})P_{k}-P_{k-1}$$
$$Q_{k}=P_{k+1}+P_k$$

If we let $P_+=P_1-qP_0$ and $P_-=P_1-q^{-1}P_0$, the solutions of these
equations can be written as follows
$$P_{k}=\frac{q^{-k}P_+-q^kP_-}{q^{-1}-q}$$
$$Q_{k}=\frac{q^{-k}P_+-q^{k+1}P_-}{1-q}$$

We can compute now the series $f_k$ by using the variables 
$z=y^{-1}=q(1+q)^{-2}$.
\begin{eqnarray*}
f_k(z)
&=&\frac{(1+q)^2}{q}\cdot\frac{P_k}{Q_k}\cr
&=&\frac{(1+q)^2}{q}\cdot\frac{1-q}{q^{-1}-q}\cdot\frac{q^{-k}P_+-q^kP_-}{q^{-k}P_+-q^{k+1}P_-}\cr
&=&(1+q)\,\frac{1-Pq^{2k}}{1-Pq^{2k+1}}\cr
\end{eqnarray*}

And finally we obtain the Jones series as
\begin{eqnarray*}
\Theta_k(q)-q
=\frac{1-q}{1+q}\,f_k(z)
=(1-q)\frac{1-Pq^{2k}}{1-Pq^{2k+1}}
\end{eqnarray*}

This proves the statement.
\end{proof}

We are now ready to compute the spectral measure for A$_n$ (on
$\mathbb{T}$).

\begin{theorem}
Let $X$ be the Coxeter-Dynkin graph $A_{n-1}$, $n \ge 1$, with $n-1$
vertices. The spectral measure of $X$ (on $\mathbb{T}$) is given by
$$d\varepsilon (u)=\alpha(u)\,d_nu$$
where $\alpha(u)=2Im(u)^2$, and $d_n$ is the uniform measure on $2n$-th 
roots of unity. 
\end{theorem}

\begin{proof}
We use Lemma 3.1 to compute the Jones series of $X_k=A_{2k+2}$ by
letting
$$M_0=\begin{pmatrix}1\end{pmatrix}{\hskip 3mm}
L_0=\begin{pmatrix}1\end{pmatrix}{\hskip 3mm}
K_0=\begin{pmatrix}\ \end{pmatrix}{\hskip 3mm}
K_1=\begin{pmatrix}2\end{pmatrix}$$
$$P_0=1$$
$$P_1=q+q^{-1}$$
$$P=q^2$$

Then
$$\frac{\Theta(q)-q}{1-q}=\frac{1-q^{2k+2}}{1-q^{2k+3}}.$$

Similarly, Lemma 3.1 gives the Jones series of $X_k=A_{2k+3}$ by
letting
$$M_0=\begin{pmatrix}1\cr 1\end{pmatrix}{\hskip 3mm}
L_0=\begin{pmatrix}1&1\cr 1&1\end{pmatrix}{\hskip 3mm}
K_0=\begin{pmatrix}1\end{pmatrix}{\hskip 3mm}
K_1=\begin{pmatrix}2&1\cr 1&1\end{pmatrix}$$
$$P_0=1+q+q^{-1}$$
$$P_1=1+q+q^{-1}+q^2+q^{-2}$$
$$P=q^{3}$$

Then
$$\frac{\Theta(q)-q}{1-q}=\frac{1-q^{2k+3}}{1-q^{2k+4}}.$$

The above two formulae give the Jones series for $A_{n-1}$ as
$$\frac{\Theta(q)-q}{1-q}=\,\frac{1-q^{n-1}}{1-q^{n}}.$$

We use the following formula, valid for $m=2nk+r$ with $r=0,1,\ldots ,2n-1$.
\begin{eqnarray*}
\int_{\mathbb{T}}\frac{u^{-m}}{1-qu}\,d_nu
&=&\frac{q^{r}}{1-q^{2n}}
\end{eqnarray*}

We can compute now the Stieltjes transform of $\varepsilon$.
\begin{eqnarray*}
2S(q)-1
&=&-1+\int_{\mathbb{T}}\frac{2-u^2-u^{-2}}{1-qu}\,d_nu\cr
&=&-1+\frac{2-q^{2n-2}-q^2}{1-q^{2n}}\cr
&=&\frac{1+q^{2n}-q^{2n-2}-q^2}{1-q^{2n}}\cr
&=&\frac{(1-q^2)(1-q^{2n-2})}{1-q^{2n}}
\end{eqnarray*}

Thus we have $2S(q)-1=\Theta(q^2)-q^2$, and we are done.
\end{proof}

We compute next the Jones series for the Coxeter-Dynkin graph of
type $D_n$, $n \ge 3$.  The graph $D_{n}$ has $n$ vertices. It consists 
of two vertices connected to one another vertex, and this vertex
in turn is connected to an $A$-tail, ending at the distinguished vertex $1$.

\begin{theorem}
Let $X$ be the Coxeter-Dynkin graph $D_{n+1}$, $n \ge 2$, with $n+1$
vertices. The spectral measure of $X$ (on $\mathbb{T}$) is given by
$$d\varepsilon (u)=\alpha(u)\,d_n^\prime u$$
where $\alpha(u)=2Im(u)^2$, and $d_n^\prime$ is the uniform measure 
on $4n$-th roots of unity of odd order. 
\end{theorem}

\begin{proof}
We use again Lemma 3.1 to compute the Jones series of $X_k=D_{2k+3}$
by letting
$$M_0=\begin{pmatrix}1&1\end{pmatrix}{\hskip 3mm}
L_0=\begin{pmatrix}2\end{pmatrix}{\hskip 3mm}
K_0=\begin{pmatrix}\ \end{pmatrix}{\hskip 3mm}
K_1=\begin{pmatrix}3\end{pmatrix}$$
$$P_0=1$$
$$P_1=-1+q+q^{-1}$$
$$P=-q$$

Thus
$$\frac{\Theta(q)-q}{1-q}=\,\frac{1+q^{2k+1}}{1+q^{2k+2}}.$$

Similarly, Lemma 3.1 allows us to compute the Jones series of 
$X_k=D_{2k+4}$ by letting
$$M_0=\begin{pmatrix}1\cr 1\cr 1\end{pmatrix}{\hskip 3mm}
L_0=\begin{pmatrix}1&1&1\cr 1&1&1\cr 1&1&1\end{pmatrix}{\hskip 3mm}
K_0=\begin{pmatrix}1&1\cr1&1\end{pmatrix}{\hskip 3mm}
K_1=\begin{pmatrix}2&1&1\cr 1&1&1\cr 1&1&1\end{pmatrix}$$
$$P_0=2+2q+2q^{-1}+q^2+q^{-2}$$
$$P_1=q+q^{-1}+2q^2+2q^{-2}+q^3+q^{-3}$$
$$P=-q^{2}$$

Thus
$$\frac{\Theta(q)-q}{1-q}=\frac{1+q^{2k+2}}{1+q^{2k+3}}.$$

From the above two formulae we deduce the Jones series of $D_{n+1}$ as
$$\frac{\Theta(q)-q}{1-q}=\,\frac{1+q^{n-1}}{1+q^{n}}.$$

We use now the following identity.
\begin{eqnarray*}
\frac{1+q^{n-1}}{1+q^{n}}
&=&\frac{(1-q^n)(1+q^{n-1})}{(1-q^{n})(1+q^n)}\cr
&=&\frac{1-q^n+q^{n-1}-q^{2n-1}}{1-q^{2n}}\cr
&=&\frac{2(1-q^{2n-1})-(1+q^n)(1-q^{n-1})}{1-q^{2n}}\cr
&=&2\,\frac{1-q^{2n-1}}{1-q^{2n}}-\frac{1-q^{n-1}}{1-q^n}
\end{eqnarray*}

Multiplying by $1-q$ and adding $q$ makes appear the Jones series 
$\Theta_m$ for $A_{m-1}$ computed in Theorem 3.1, and we obtain
$$\Theta(q)=2\Theta_{2n}(q)-\Theta_n(q).$$

We compute now the Stieltjes transform of $\varepsilon$. Denote
by $S_m$ the Stieltjes transform of the measure for $A_{m-1}$.
Then
\begin{eqnarray*}
2S(q)-1
&=&2(2S_{2n}(q)-S_n(q))-1\cr
&=&2(S_{2n}(q)-1)-(2S_n(q)-1)\cr
&=&2(\Theta_{2n}(q^2)-q^2)-(\Theta_n(q^2)-q^2)\cr
&=&(2\Theta_{2n}(q^2)-\Theta_n(q^2))-q^2
\end{eqnarray*}

Thus we have $2S(q)-1=\Theta(q^2)-q^2$ as claimed.
\end{proof}

\begin{theorem}
For $X=A_\infty$ the spectral measure is given by
$$d\varepsilon (u)=\alpha(u)\,du$$
where $\alpha(u)=2Im(u)^2$, and $du$ is the uniform measure on the 
unit circle.
\end{theorem}

\begin{proof}
This follows from Theorems 3.1 or 3.2 with $n\to\infty$, or from a 
direct loop count.
\end{proof}

\section{Graphs with fork tails}

Consider a sequence of graphs $X_k$ obtained by adding $D_{2k+2}$ tails 
to a given graph (compare with section 3). We let
$$X_0=\begin{matrix}1\cr 2\end{matrix}>\alpha:::\Gamma$$
$$X_1=\begin{matrix}1\cr 2\end{matrix}>\alpha-3-\beta:::\Gamma$$
$$X_2=\begin{matrix}1\cr 2\end{matrix}>\alpha-3-\beta-4-\gamma:::\Gamma$$

As before we let $L_0$ be the matrix appearing on the top left of the 
adjacency matrix of $X_0$, i.e.
$$\Delta_0^2=\begin{pmatrix}L_0&0\cr 0&N_0\end{pmatrix}.$$

\begin{lemma}
The Jones series of the graphs $X_{k}$ is given by
$$(\Theta(q)-q)(1+q)=\frac{1-Pq^{2k+1}}{1+Pq^{2k}}$$
where $P$ is defined by 
$$P=\frac{P_1-q^{-1}P_0}{P_1-qP_0}$$
where $P_i=P_i(y)=\det (y-J_i)$, $i=0$, $1$, with $y=2+q+q^{-1}$, and with 
$J_0$, $J_1$ being the following matrices:

-- $J_0$ is obtained from $L_0$ by deleting the first two rows and columns.

-- $J_1$ is obtained from $L_0$ by deleting the first row and column, then 
adding $1$ to the first entry.
\end{lemma}

\begin{proof}
The matrix $M_k$ associated to the adjacency matrix of $X_k$ is described by a 
column vector $w$, and a matrix $M$, as follows:
$$M_0=\begin{pmatrix}1&\cr 1&\cr w&M\end{pmatrix}\hskip 5mm 
M_1=\begin{pmatrix}1&&\cr1&&\cr1&1&\cr&w&M\end{pmatrix}\hskip 5mm 
M_2=\begin{pmatrix}1&&&\cr1&&&\cr1&1&&\cr&1&1&\cr&&w&M\end{pmatrix}$$

The corresponding matrices $L_k$ make appear the matrix $N=ww^t+MM^t$
and are given by
$$L_0=\begin{pmatrix}1&\cr 1&\cr w&M\end{pmatrix}
\begin{pmatrix}1&1&w^t\cr &&M^t\end{pmatrix}=
\begin{pmatrix}1&1&w^t\cr1&1&w^t\cr w&w&N\end{pmatrix}$$
$$L_1=
\begin{pmatrix}1&&\cr1&&\cr1&1&\cr&w&M\end{pmatrix}
\begin{pmatrix}1&1&1&\cr&&1&w^t\cr&&&M^t\end{pmatrix}=
\begin{pmatrix}1&1&1&\cr1&1&1&\cr1&1&2&w^t\cr&&w&N\end{pmatrix}$$
$$L_2=
\begin{pmatrix}1&&&\cr1&&&\cr1&1&&\cr&1&1&\cr&&w&M\end{pmatrix}
\begin{pmatrix}1&1&1&&\cr&&1&1&\cr&&&1&w^t\cr&&&&M^t\end{pmatrix}=
\begin{pmatrix}1&1&1&&\cr1&1&1&&\cr1&1&2&1&\cr&&1&2&w^t\cr&&&w&N\end{pmatrix}$$

It is now clear what the form of the matrices $M_k$ and $L_k$ is for 
general $k$. Consider the matrix $K_k$, obtained from $L_k$ by deleting 
the first row and column, that is
$$K_0=\begin{pmatrix}1&w^t\cr w&N\end{pmatrix}\hskip 5mm K_1=
\begin{pmatrix}1&1&\cr1&2&w^t\cr&w&N\end{pmatrix}\hskip 5mm K_2=
\begin{pmatrix}1&1&&\cr1&2&1&\cr&1&2&w^t\cr&&w&N\end{pmatrix}$$
and similarly for general $k$.
Consider also the matrix $J_k$, obtained from $L_k$ by deleting the 
first two rows and columns, i.e.
$$J_0=\begin{pmatrix}N\end{pmatrix}\hskip 5mm 
J_1= \begin{pmatrix}2&w^t\cr w&N\end{pmatrix}\hskip 5mm 
J_2= \begin{pmatrix}2&1&\cr1&2&w^t\cr&w&N\end{pmatrix}$$

We have then the following formula for the Poincar\'e series of $X_k$, 
with $y=z^{-1}$.
\begin{eqnarray*}
f_k(z)
=\int\frac{1}{1-zL_k}
=\frac{\det (1-zK_k)}{\det (1-zL_{k})}
=y\cdot\frac{\det (y-K_{k})}{\det (y-L_{k})}\cr
\end{eqnarray*}

The characteristic polynomials $P_k=\det (y-J_k)$, $Q_k=\det (y-K_k)$ and 
$R_k=\det(y-L_k)$ satisfy the following relations, obtained by developing 
the determinant at top left:
$$P_{k+1}=(y-2)P_{k}-P_{k-1}$$
$$Q_{k}=(y-1)P_{k}-P_{k-1}$$
$$R_k=(y-1)Q_k-P_{k}-(y+1)P_{k-1}$$

The solutions of these equations can be written in terms of 
$P_+=P_1-qP_0$ and $P_-=P_1-q^{-1}P_0$, namely 
$$P_{k}=\frac{q^{-k}P_+-q^kP_-}{q^{-1}-q}$$
$$Q_{k}=\frac{q^{-k}P_+-q^{k+1}P_-}{1-q}$$
$$R_k=\frac{q^{-k}P_++q^kP_-}{q/(1+q)^2}$$

We can compute now $f_k$ by using the variables $z=y^{-1}=q(1+q)^{-2}$.
\begin{eqnarray*}
f_k(z)
&=&\frac{(1+q)^2}{q}\cdot\frac{Q_k}{R_k}\cr
&=&\frac{(1+q)^2}{q}\cdot\frac{q/(1+q)^2}{1-q}\cdot\frac{q^{-k}P_+-q^{k+1}P_-}{q^{-k}P_++q^{k}P_-}\cr
&=&\frac{1}{1-q}\cdot\frac{1-Pq^{2k+1}}{1+Pq^{2k}}\cr
\end{eqnarray*}

And finally we obtain the Jones series for $D_{n+1}$ as
\begin{eqnarray*}
\Theta_k(q)-q
=\frac{1-q}{1+q}\,f_k(z)
=\frac{1}{1+q}\cdot\frac{1-Pq^{2k+1}}{1+Pq^{2k}}
\end{eqnarray*}

This proves the claim.
\end{proof}

We proceed now with the computation of the spectral measure for the
extended Coxeter-Dynkin graph $D_n^{(1)}$.
The graph $D_n^{(1)}$, $n\geq 4$, has $n+1$ vertices. It consists of 
the distinguished vertex $1$, connected to a triple point, connected in 
turn to another vertex and to an $A$-tail ending at another triple point.
This triple point is connected to two other vertices.

\begin{theorem}
Let $X$ be the extended Coxeter-Dynkin graph of type $D_{n+2}^{(1)}$,
$n \ge 2$ ($n+3$ vertices). The spectral measure of $X$
(on $\mathbb{T}$) is given by
$$d\varepsilon(u)=\frac{\delta_i+\delta_{-i}}{4}+\frac{d_nu}{2},$$
where $d_nu$ is the uniform measure on $2n$-th roots of unity and
$\delta_w$ is the Dirac measure at $w \in \mathbb{T}$.
\end{theorem}

\begin{proof}
We use Lemma 4.1 to compute the Jones series of $X_k=D_{2k+4}^{(1)}$
by letting
$$M_0=\begin{pmatrix}1\cr 1\cr 1\cr 1\end{pmatrix}{\hskip 3mm}
L_0=\begin{pmatrix}1&1&1&1\cr1&1&1&1\cr1&1&1&1\cr1&1&1&1\end{pmatrix}
{\hskip 3mm}
J_0=\begin{pmatrix}1&1\cr1&1\end{pmatrix}{\hskip 3mm}
J_1=\begin{pmatrix}2&1&1\cr1&1&1\cr1&1&1\end{pmatrix}$$
$$P_0=2+2q+2q^{-1}+q^2+q^{-2}$$
$$P_1=q+q^{-1}+2q^2+2q^{-2}+q^3+q^{-3}$$
$$P=-q^{2}$$

We obtain
$$(\Theta(q)-q)(1+q)=\frac{1+q^{2k+3}}{1-q^{2k+2}}.$$

Similarly, Lemma 4.1 gives the Jones series of $X_k=D_{2k+5}^{(1)}$
by letting
$$M_0=\begin{pmatrix}1&0&0\cr 1&0&0\cr 1&1&1\end{pmatrix}{\hskip 3mm}
L_0=\begin{pmatrix}1&1&1\cr 1&1&1\cr 1&1&3\end{pmatrix}{\hskip 3mm}
K_0=\begin{pmatrix}3\end{pmatrix}{\hskip 3mm}
K_1=\begin{pmatrix}2&1\cr 1&3\end{pmatrix}$$
$$P_0=-1+q+q^{-1}$$
$$P_1=1-q-q^{-1}+q^2+q^{-2}$$
$$P=-q^3$$

We obtain
$$(\Theta(q)-q)(1+q)=\frac{1+q^{2k+4}}{1-q^{2k+3}}.$$

From the above two formulae we deduce that Jones series for 
$D_{n+2}^{(1)}$ satisfies
$$(\Theta(q)-q)(1+q)=\frac{1+q^{n+1}}{1-q^{n}}.$$

Hence we get the following explicit formula for the Jones series
for $D_{n+2}^{(1)}$:
\begin{eqnarray*}
\Theta (q)-q
&=&\frac{1}{1+q}\cdot \frac{1-q^n+q^n+q^{n+1}}{1-q^n}\cr
&=&\frac{1}{1+q}\left( 1+\frac{q^n(1+q)}{1-q^n}\right)\cr
&=&\frac{1}{1+q}+\frac{q^n}{1-q^n}\cr
&=&\frac{1}{1+q}+\frac{1}{1-q^n}-1
\end{eqnarray*}

We compute now the Stieltjes transform of the associated spectral
measure $\varepsilon$.
\begin{eqnarray*}
2S(q)
&=&\frac{1}{2}\left(\frac{1}{1-qi}+\frac{1}{1+qi}\right) +
\int_{\mathbb{T}}\frac{1}{1-qu}\,d_nu\cr
&=&\frac{1}{1+q^2}+\frac{1}{1-q^{2n}}
\end{eqnarray*}

Thus we have $2S(q)=\Theta (q^2)-q^2+1$, and we are done.
\end{proof}

The spectral measure for the infinite Coxeter-Dynkin graph $D_\infty$
follows now. Recall that the graph $D_\infty$ has a triple point, 
connected to the distinguished vertex $1$ and to another vertex, and to 
an A$_\infty$-tail.

\begin{theorem}
The spectral measure (on $\mathbb{T}$) of the Coxeter-Dynkin graph 
$D_\infty$ is given by
$$d\varepsilon(u)=\frac{\delta_i+\delta_{-i}}{4}+\frac{du}{2},$$
where $du$ is the uniform measure on the unit circle. 
\end{theorem}

\begin{proof}
This follows from Theorem 4.1 by letting $n\to\infty$.
\end{proof}

\section{Exceptional graphs}

In this chapter we computer the spectral measures for the
Coxeter-Dynkin graphs of type $E$. These graphs are
\begin{eqnarray*}
E_6&=&F(2,1,2)\cr
E_7&=&F(2,1,3)\cr
E_8&=&F(2,1,4)\cr
E_6^{(1)}&=&F(2,2,2)\cr
E_7^{(1)}&=&F(3,1,3)\cr 
E_8^{(1)}&=&F(2,1,5)
\end{eqnarray*}
Here we denote by $F(a,b,c)$ the graph with $a+b+c+1$ vertices, 
consisting of a triple point which is connected to an $A_a$ tail, 
to an $A_b$ tail, and to an $A_c$ tail ending at the distinguished 
vertex $1$.
All of these graphs, with the exception of $E_7$, appear as principal
graphs of subfactors with Jones index $\le 4$ (see for instance 
\cite{ek}, \cite{ghj}).

We use the following notation. Assume that $\varepsilon$ is a probability 
measure on the unit circle which is even, in the sense that all its odd 
moments are $0$. The Stieltjes transform $S(q)$ of $\varepsilon$ is then
a series in $q^2$, and the following definition makes sense.

\begin{definition}
The $T$ series of an even measure $\varepsilon$ is given by
$$T(q)=\frac{2S(q^{1/2})-1}{1-q}$$
where $S$ is the Stieltjes transform of $\varepsilon$.
\end{definition}

It follows from Propositions 1.1 and 1.2 that the spectral measure of 
a graph is even, and that we have the formula
$$T(q)=\frac{\Theta(q)-q}{1-q}$$
where $\Theta$ is the Jones series. 

In this section we compute the $T$ series of exceptional graphs of type E.

\begin{theorem}
The $T$ series of the Coxeter-Dynkin graphs $E_6$, $E_7$ and $E_8$ are 
given by the following formulae:
\begin{eqnarray*}
T_6(q)&=&\frac{(1-q^6)(1-q^8)}{(1-q^3)(1-q^{12})}\cr
T_7(q)&=&\frac{(1-q^9)(1-q^{12})}{(1-q^4)(1-q^{18})}\cr
T_8(q)&=&\frac{(1-q^{10})(1-q^{15})(1-q^{18})}{(1-q^5)(1-q^9)(1-q^{30})}
\end{eqnarray*}
\end{theorem}

The proof of these results uses techniques from section 3. We combine
the proof with the proof of next result.

\begin{theorem}
The $T$ series of the extended Coxeter-Dynkin graphs 
$E_6^{(1)}$, $E_7^{(1)}$ and $E_8^{(1)}$ are given by the following formulae:
\begin{eqnarray*}
T_6^{(1)}(q)&=&\frac{1-q^{12}}{(1-q^3)(1-q^4)(1-q^6)}\cr T_7^{(1)}(q)&=&\frac{1-q^{18}}{(1-q^4)(1-q^6)(1-q^9)}\cr T_8^{(1)}(q)&=&\frac{1-q^{30}}{(1-q^6)(1-q^{10})(1-q^{15})}
\end{eqnarray*}
\end{theorem}

\begin{proof}
We compute first the $T$ series of $E_6,E_7,E_8,E_8^{(1)}$, which are all 
of the form $F(2,1,n)$. We use Lemma 3.1 to compute the $T$ series 
of $X_k=F(2,1,2k)$ by letting
$$M_0=\begin{pmatrix}1&1\cr 0&1\end{pmatrix}{\hskip 3mm}
L_0=\begin{pmatrix}2&1\cr 1&1\end{pmatrix}{\hskip 3mm}
K_0=\begin{pmatrix}1\end{pmatrix}{\hskip 3mm}
K_1=\begin{pmatrix}3&1\cr 1&1\end{pmatrix}$$
$$P_0=1+q+q^{-1}$$
$$P_1=q^2+q^{-2}$$
$$P=-q\,\frac{1+q-q^3}{1-q^2-q^3}$$

We obtain
$$T(q)=\frac{(1-q^2-q^{3})+q^{2k+1}(1+q-q^3)}{(1-q^2-q^{3})+
q^{2k+2}(1+q-q^3)}.$$

Similarly, Lemma 3.1 gives the $T$ series of $X_k=F(2,1,2k+1)$ by letting
$$M_0=\begin{pmatrix}1&0\cr 1&0\cr 1&1\end{pmatrix}{\hskip 3mm}
L_0=\begin{pmatrix}1&1&1\cr 1&1&1\cr 1&1&2\end{pmatrix}{\hskip 3mm}
K_0=\begin{pmatrix}1&1\cr 1&2\end{pmatrix}{\hskip 3mm}
K_1=\begin{pmatrix}2&1&1\cr 1&1&1\cr 1&1&2\end{pmatrix}$$
$$P_0=1+q+q^{-1}+q^2+q^{-2}$$
$$P_1=-1+q^2+q^{-2}+q^3+q^{-3}$$
$$P=-q^2\,\frac{1+q-q^3}{1-q^2-q^3}$$

We obtain
$$T(q)=\frac{(1-q^2-q^{3})+q^{2k+2}(1+q-q^3)}{(1-q^2-q^{3})+
q^{2k+3}(1+q-q^3)}.$$

To simplify this expression we introduce the following polynomials:
$$Q_n=(1-q^2-q^{3})+q^{n}(1+q-q^3)$$

From the above formulae we get then the $T$ series of $F(2,1,n)$ in 
terms of these polynomials as
$$T(q)=\frac{Q_{n+1}}{Q_{n+2}}.$$

The polynomials $Q_k$ needed for $E_6$, $E_7$, $E_8$,
and $E_8^{(1)}$ are all cyclotomic:
\begin{eqnarray*}
Q_3&=&\frac{(1-q^2)(1-q^8)}{1-q^4}\\
Q_4&=&\frac{(1-q^2)(1-q^3)(1-q^{12})}{(1-q^4)(1-q^6)}\\
Q_5&=&\frac{(1-q^2)(1-q^3)(1-q^{18})}{(1-q^6)(1-q^9)}\\
Q_6&=&\frac{(1-q^2)(1-q^3)(1-q^5)(1-q^{30})}{(1-q^6)(1-q^{10})(1-q^{15})}\\
Q_7&=&(1-q^2)(1-q^3)(1-q^5)
\end{eqnarray*}

The formulae for the $T$ series of $E_6$, $E_7$, $E_8$ and $E_8^{(1)}$
follow now.

We use again Lemma 3.1 to compute the $T$ series of $X_k=F(2,2,2k)$ by
letting
$$M_0=\begin{pmatrix}1&1\cr 1&0\cr 0&1\end{pmatrix}{\hskip 3mm}
L_0=\begin{pmatrix}2&1&1\cr 1&1&0\cr 1&0&1\end{pmatrix}{\hskip 3mm}
K_0=\begin{pmatrix}1&0\cr 0&1\end{pmatrix}{\hskip 3mm}
K_1=\begin{pmatrix}3&1&1\cr 1&1&0\cr 1&0&1\end{pmatrix}$$
$$P_0=3+2q+2q^{-1}+q^2+q^{-2}$$
$$P_1=-1+q^2+q^{-2}+q^3+q^{-3}$$
$$P=-q\,\frac{1+q-q^2}{1-q-q^2}$$

We obtain
$$T(q)=\frac{(1-q-q^2)+q^{2k+1}(1+q-q^2)}{(1-q-q^2)+q^{2k+2}(1+q-q^2)}.$$

Factorising the numerator $N$ and denominator $D$ for $k=1$, we get
$$N=1-q-q^2+q^3+q^4-q^5=(1-q)(1-q^2+q^4),$$
$$D=1-q-q^2+q^4+q^5-q^6=(1-q)(1-q^2)(1-q^3).$$

Rewrite these expressions as
\begin{eqnarray*}
N&=&\frac{(1-q)(1-q^2)(1-q^{12})}{(1-q^4)(1-q^6)}\cr
D&=&(1-q)(1-q^2)(1-q^3)
\end{eqnarray*}
and the formula for the $T$ series of $E_6^{(1)}$ follows.

We use next Lemma 3.1 to compute the $T$ series of $X_k=F(3,1,2k+1)$
by letting
$$M_0=\begin{pmatrix}1&0\cr 1&0\cr 1&1\cr 0&1\end{pmatrix}{\hskip 3mm}
L_0=\begin{pmatrix}1&1&1&0\cr 1&1&1&0\cr 1&1&2&1\cr 0&0&1&1\end{pmatrix}{\hskip 3mm}
K_0=\begin{pmatrix}1&1&0\cr 1&2&1\cr 0&1&1\end{pmatrix}{\hskip 3mm}
K_1=\begin{pmatrix}2&1&1&0\cr 1&1&1&0\cr 1&1&2&1\cr 0&0&1&1\end{pmatrix}$$
$$P_0=2+2q+2q^{-1}+2q^2+2q^{-2}+q^3+q^{-3}$$
$$P_1=-2-q-q^{-1}+q^2+q^{-2}+2q^3+2q^{-3}+q^4+q^{-4}$$
$$P=-q^2\,\frac{1+q^2-q^3}{1-q-q^3}$$

We obtain
$$T(q)=\frac{(1-q-q^3)+q^{2k+2}(1+q^2-q^3)}{(1-q-q^3)+q^{2k+3}(1+q^2-q^3)}.$$

Factorising the numerator $N$ and denominator $D$ for $k=1$ gives
$$N=1-q-q^3+q^4+q^6-q^7=(1-q)(1-q^3+q^6),$$
$$D=1-q-q^3+q^5+q^7-q^8=(1-q)(1-q^3)(1-q^4).$$

Rewrite these expressions as
\begin{eqnarray*}
N&=&\frac{(1-q)(1-q^3)(1-q^{18})}{(1-q^6)(1-q^9)}\cr
D&=&(1-q)(1-q^3)(1-q^4)
\end{eqnarray*}
and the formula for the $T$ series of $E_7^{(1)}$ follows.
\end{proof}

\section{Spectral measures for exceptional graphs}

We compute in this section the spectral measures of the exceptional 
graphs $E_6$ and $E_{6,7,8}^{(1)}$. We will express the $T$ series 
computed in the previous 
section as linear combinations of elementary $T$ series computed
in the next two Lemmas. The notation for the measures is as in sections
3 and 4.

\begin{lemma}
The $T$ series of the measures $\alpha\,d_n$, $\alpha\,d_n'$, $d_n$, 
$d_n'$ are given by the following identities:
\begin{eqnarray*}
T(\alpha\,d_n)&=&\frac{1-q^{n-1}}{1-q^n}\cr
T(\alpha\,d_n')&=&\frac{1+q^{n-1}}{1+q^n}\cr
T(d_n)&=&\frac{1}{1-q}\cdot\frac{1+q^n}{1-q^n}\cr
T(d_n')&=&\frac{1}{1-q}\cdot\frac{1-q^n}{1+q^n}
\end{eqnarray*}
\end{lemma}

\begin{proof}
The first two identities appear in proof of Theorems 3.1 and 3.2. For 
the third one we use that
$$\int_{\mathbb{T}} u^kd_nu=(2n|k),$$
where $(2n|k)$ is defined to be $1$ when $2n$ divides $k$, and $0$ 
otherwise. This gives the following identity for the Stieltjes transform:
$$S(q)=\sum_{s=0}^\infty q^{2ns}=\frac{1}{1-q^{2n}}.$$

We can then compute the $T$ series as in Definition 5.1 by
\begin{eqnarray*}
T(d_n)
&=&\frac{2S(q^{1/2})-1}{1-q}\cr
&=&\frac{1}{1-q}\left(\frac{2}{1-q^n}-1\right)\cr
&=&\frac{1}{1-q}\cdot\frac{1+q^n}{1-q^n}
\end{eqnarray*}

For the fourth identity in the lemma we use that $d_n'=2d_{2n}-d_n$.
Hence
\begin{eqnarray*}
T(d_n')
&=&2T(d_{2n})-T(d_n)\cr
&=&\frac{1}{1-q}\left(2\cdot\frac{1+q^{2n}}{1-q^{2n}}-\frac{1+q^n}{1-q^n}\right)\cr
&=&\frac{1}{1-q}\cdot\frac{1-2q^n+q^{2n}}{1-q^{2n}}\cr
&=&\frac{1}{1-q}\cdot\frac{(1-q^n)^2}{1-q^{2n}}
\end{eqnarray*}

After simplification we obtain the formula in the statement.
\end{proof}

\begin{theorem}
The spectral measures of $E_{6,7,8}^{(1)}$ (on $\mathbb{T}$) are given 
by
\begin{eqnarray*}
\varepsilon_6^{(1)}&=&\alpha\,d_3+(d_2-d_3)/2\cr
\varepsilon_7^{(1)}&=&\alpha\,d_4+(d_3-d_4)/2\cr
\varepsilon_8^{(1)}&=&\alpha\,d_6+(d_5-d_6)/2
\end{eqnarray*}
where $\alpha(u)=2Im(u)^2$, and $d_nu$ is the uniform measure on 
$2n$-th roots of unity.
\end{theorem}

\begin{proof}
The $T$ series in Theorem 5.2 can be written as
\begin{eqnarray*}
T_6^{(1)}&=&\frac{1+q^6}{(1-q^3)(1-q^4)}\cr T_7^{(1)}&=&\frac{1+q^9}{(1-q^4)(1-q^6)}\cr T_8^{(1)}&=&\frac{1+q^{15}}{(1-q^6)(1-q^{10})}
\end{eqnarray*}

Factoring by $1+q^2$, $1+q^3$ resp. $1+q^5$ gives
\begin{eqnarray*}
T_6^{(1)}&=&\frac{1-q^2+q^4}{(1-q^2)(1-q^3)}\cr T_7^{(1)}&=&\frac{1-q^3+q^6}{(1-q^3)(1-q^4)}\cr T_8^{(1)}&=&\frac{1-q^5+q^{10}}{(1-q^5)(1-q^{6})}
\end{eqnarray*}

We get then the following formula, with $k=2,3,5$ corresponding to 
$n=6,7,8$:
$$T_{n}^{(1)}
=\frac{1-q^{k}+q^{2k}}{(1-q^{k})(1-q^{k+1})}.$$

We can rewrite this series in the following way:
\begin{eqnarray*}
T_{n}^{(1)}
&=&\frac{1-2q^{k}+q^{k}}{(1-q^{k})(1-q^{k+1})}+\frac{q^{k}}{(1-q^{k})(1-q^{k+1})}\cr
&=&\frac{1-q^{k}}{1-q^{k+1}}+\frac{1}{1-q}\cdot\frac{q^{k}-q^{k+1}}{(1-q^{k})(1-q^{k+1})}\cr
&=&\frac{1-q^{k}}{1-q^{k+1}}+\frac{1}{1-q}\left( \frac{1}{1-q^{k}}-\frac{1}{1-q^{k+1}}\right)
\end{eqnarray*}

We can then write $T_{n}^{(1)}$ as a linear combination of the 
elementary $T$ series from Lemma 6.1 as follows
$$T_n^{(1)}=T(\alpha\,d_{k+1})+(T(d_{k})-T(d_{k+1}))/2.$$

By using linearity of the Stieltjes transform, hence of the $T$ series, 
we get the formulae in the statement of the theorem.
\end{proof}

\begin{theorem}
The spectral measure of $E_{6}$ (on $\mathbb{T}$) is given by
$$\varepsilon_6=\alpha\,d_{12}+(d_{12}-d_6-d_4+d_3)/2$$
where $\alpha(u)=2Im(u)^2$, and $d_n$ is the uniform measure on 
$2n$-th roots of unity.
\end{theorem}

\begin{proof}
The $T$ series of $E_6$ can be written as
\begin{eqnarray*}
T_6
&=&\frac{(1+q^3)(1-q^8)}{1-q^{12}}\cr
&=&\frac{1-q^{11}}{1-q^{12}}+\frac{q^3-q^8}{1-q^{12}}.
\end{eqnarray*}

Note that we have the following identity:
\begin{eqnarray*}
\frac{q^3-q^8}{1-q^{12}}
&=&\frac{1}{1-q}\cdot\frac{q^3-q^4-q^8+q^9}{1-q^{12}}\cr
&=&\frac{1}{1-q}\left(\frac{1}{1-q^{12}}-\frac{1+q^6}{1-q^{12}}-\frac{1+q^4+q^8}{1-q^{12}}+\frac{1+q^3+q^6+q^9}{1-q^{12}}\right)\cr
&=&\frac{1}{1-q}\left(\frac{1}{1-q^{12}}-\frac{1}{1-q^6}-\frac{1}{1-q^4}+
\frac{1}{1-q^3}\right).
\end{eqnarray*}

It follows now easily that $T_6$ can be written as a linear 
combination of the elementary $T$ series in Lemma 6.1, namely
$$T_6=T(\alpha\,d_{12})+(T(d_{12})-T(d_6)-T(d_4)+T(d_3))/2.$$

This gives the formula for the spectral measure $\varepsilon_6$ in 
the statement of the theorem.
\end{proof}

\section{Exceptional measures: $E_7$ and $E_8$}

Note that all spectral measures of the finite ADE graphs computed 
so far are linear combinations of measures of type $d_n$ and 
$\alpha\, d_n$ (observe that we have $d_n'=2d_{2n}-d_n$).  

\begin{definition}
A discrete measure supported by roots of unity is called {\rm cyclotomic} 
if it is a linear combination of measures of type $d_n$, $n\geq 1$, and 
$\alpha\, d_n$, $n\geq 2$.
\end{definition}

Note that we require $n\geq 2$ for the measure $\alpha\, d_n$. This is 
simply because $\alpha\, d_1$ is the null measure.

\begin{theorem}
The spectral measures of $E_7$, $E_8$ (on $\mathbb{T}$) are not 
cyclotomic.
\end{theorem}

\begin{proof}
From Theorem 5.1 we obtain the $T$ series of $E_7$ as
$$T_7=\frac{(1-q^9)(1+q^4+q^8)}{1-q^{18}}.$$

This shows that the corresponding spectral measure $\varepsilon_7$ is 
supported by $36$-th roots of unity. Assume now that 
$\varepsilon_7$ is cyclotomic. Then
$$\varepsilon_7\in {\rm span}\{ d_n,\,\alpha\,d_m\mid n\geq 1,\, 
m\geq 2,\, n,m|18\}.$$

Using the linearity of the $T$ transform (with respect to the measure)
this means
$$T_7\in {\rm span}\{ T(d_n),\,T(\alpha\,d_m)\mid n\geq 1,\, 
m\geq 2,\, n,m|18\}.$$

We multiply the
relevant $T$ series by $(1-q)(1-q^{18})$, i.e. we consider the 
following degree $18$ polynomials, where $n,m$ are as above:
\begin{eqnarray*}
P_n&=&(1-q)(1-q^{18})\,T(d_n)\cr
Q_m&=&(1-q)(1-q^{18})\,T(\alpha\,d_m)\cr
R_7&=&(1-q)(1-q^{18})\,T_7
\end{eqnarray*}

The assumption that $\varepsilon_7$ is cyclotomic becomes then
$$R_7\in {\rm span}\{ P_n,\,Q_m\mid n\geq 1,\, m\geq 2,\, n,m|18\}.$$

The above formula of $T$ and Lemma 6.1 lead to the following expessions:
\begin{eqnarray*}
P_n&=&(1+q^{n})\cdot\frac{1-q^{18}}{1-q^n}\cr
Q_m&=&(1-q)(1-q^{m-1})\cdot\frac{1-q^{18}}{1-q^m}\cr
R_7&=&(1-q)(1-q^{9})(1+q^4+q^8)
\end{eqnarray*}

The coefficients $c_k$ of each of these polynomials satisfy 
$c_k=c_{18-k}$, so in order to solve the system of 
linear equations resulting from our assumption, we can restrict 
attention to coefficients $c_k$ 
with $k=0,1,\ldots ,9$. Thus we have $10$ equations, and the unknowns 
are the coefficients of $P_n,Q_m$ with $n\geq 1$, $m\geq 2$ and $n,m|18$.

The matrix of the system is given below. The rows correspond to 
the polynomials appearing, and the columns correspond to coefficients 
of $q^k$, with $k=0,1,\ldots ,9$:
$$\begin{array}{ccccccccccc}
        \ \ \ \ \ \ \ \  &c_0&c_1&c_2&c_3&c_4&c_5&c_6&c_7&c_8&c_9\\
        &\\
 
 P_1   &1& 2& 2& 2& 2& 2& 2& 2& 2& 2\\
 P_2   &1&  & 2&  & 2&  & 2&  & 2&  \\
 P_3   &1&  &  & 2&  &  & 2&  &  & 2\\
 P_6   &1&  &  &  &  &  & 2&  &  &  \\
 P_9   &1&  &  &  &  &  &  &  &  & 2\\
 P_{18}&1&  &  &  &  &  &  &  &  &  \\
 Q_2   &1&-2& 2&-2& 2&-2& 2&-2& 2&-2\\
 Q_3   &1&-1&-1& 2&-1&-1& 2&-1&-1& 2\\
 Q_6   &1&-1&  &  &  &-1& 2&-1&  &  \\
 Q_9   &1&-1&  &  &  &  &  &  &-1& 2\\
 Q_{18}&1&-1&  &  &  &  &  &  &  &  \\
% &\\
 R_7     &1&-1&  &  & 1&-1&  &  & 1&-2\\
 \end{array}$$
 
Comparing the $c_2$ and $c_4$ columns shows that this system of 
equations has no solution. This contradicts our assumption 
that $\varepsilon_7$ is cyclotomic.

The same method applies to $E_8$. From Theorem 5.1 we get the 
$T$ series for $E_8$ as
$$T_8=\frac{(1+q^5)(1+q^9)(1-q^{15})}{1-q^{30}}.$$

Thus $\varepsilon_8$ is supported by $60$-th roots of unity. Now by 
using degree $30$ polynomials $P_n,Q_m$ and $R_8$ defined as above, 
we get again a linear system of equations with the following matrix 
of coefficients:
\begin{small}
$$\begin{array}{ccccccccccccccccc}
        \ \ \ \ \ \ \ \  &c_0&c_1&c_2&c_3&c_4&c_5&c_6&c_7&c_8&c_9&c_{10}&c_{11}&c_{12}&c_{13}&c_{14}&c_{15}\\
        &\\
 
 P_1   &1& 2& 2& 2& 2& 2& 2& 2& 2& 2&2&2&2&2&2&2\\
 P_2   &1&  & 2&  & 2&  & 2&  & 2&  &2&&2&&2&\\
 P_3   &1&  &  & 2&  &  & 2&  &  & 2&&&2&&&2\\
 P_5   &1&  &  &  &  &  2& &  &  &  &2&&&&&2\\
 P_6   &1&  &  &  &  &  & 2&  &  &  &&&2&&&\\
 P_{10}   &1&  &  &  &  &  &  &  &  & &2&&&&&\\
 P_{15}&1&  &  &  &  &  &  &  &  &  &&&&&&2\\
 P_{30}&1&  &  &  &  &  &  &  &  &  \\
 Q_2   &1&-2& 2&-2& 2&-2& 2&-2& 2&-2&2&-2&2&-2&2&-2\\
 Q_3   &1&-1&-1& 2&-1&-1& 2&-1&-1& 2&-1&-1&2&-1&-1&2\\
 Q_5   &1&-1&  &  &  -1& 2&-1&  &&-1&2&-1&&&-1&2  \\
 Q_6   &1&-1&  &  &  &-1& 2&-1&  &  &&-1&2&-1&&\\
 Q_{10}   &1&-1&  &  &  &  &  &  &&-1& 2&-1\\
 Q_{15}&1&-1&  &  &  &  &  &  &  &&&&&&-1&2  \\
 Q_{30}&1&-1&  &  &  &  &  &  &  &  \\
% &\\
 R_8     &1&-1&  &  && 1&-1&  &  & 1&-1&&&&1&-2\\
 \end{array}$$
 \end{small}
 
Assume that $\varepsilon_8$ is cyclotomic. This means that $R_8$ 
appears as linear combination of $P_n$, $Q_m$. Now, comparing 
the $c_2$ and $c_4$ columns shows that the coefficient of $Q_5$ must 
be zero, and comparing the $c_6$ and $c_{12}$ columns shows that 
the coefficient of $Q_5$ must be non-zero. Thus our assumption 
that $\varepsilon_8$ is cyclotomic is wrong.
\end{proof}

\end{document}